\documentclass[11pt]{amsart}
\usepackage{amsmath,amssymb,txfonts}
\usepackage{amssymb}
\usepackage{amsmath}
\usepackage{mathrsfs}
\usepackage[shortlabels]{enumitem}

\usepackage{color}
\usepackage{extarrows}
\newtheorem{theorem}{Theorem}[section]
\newtheorem{lemma}[theorem]{Lemma}
\newtheorem{proposition}[theorem]{Proposition}
\newtheorem{corollary}[theorem]{Corollary}
\theoremstyle{definition}
\newtheorem{definition}[theorem]{Definition}

\theoremstyle{remark}
\newtheorem{remark}[theorem]{Remark}
\numberwithin{equation}{section}

\begin{document}

\title[Fractional Sobolev Inequalities]{Strengthened   Fractional Sobolev Type Inequalities in Besov Spaces}

\author[P. Li]{Pengtao Li}
\address[P. Li]{Department of Mathematics, Qingdao University, Qingdao, Shandong 266071, China}
\email{ptli@qdu.edu.cn}

\author[R. Hu]{Rui Hu}
\address[R. Hu]{Department of Mathematics and Statistics,  MacEwan University   Edmonton, Alberta T5J2P2 Canada}
\email{hur3@macewan.ca}

\author[Z. Zhai]{Zhichun Zhai}
\address[Z. Zhai]{Department of Mathematics and Statistics,  MacEwan University,   Edmonton, Alberta T5J2P2 Canada
}
\email{zhaiz2@macewan.ca}

\thanks{Corresponding author:  zhaiz2@macewan.ca}
\thanks{Project supported:  Pengtao Li was  supported by National Natural Science Foundation of China (No. 11871293, No. 12071272) and Shandong Natural Science Foundation of China (No. ZR2017JL008).}

\date{}

\dedicatory{}


\begin{abstract} The purpose of this article is twofold. The first is to strengthen fractional Sobolev type inequalities in Besov spaces via the classical Lorentz space. In doing so, we show that the Sobolev inequality in Besov spaces is equivalent to the fractional Hardy inequality and the iso-capacitary type inequality. Secondly, we will 
strengthen fractional Sobolev type inequalities in Besov spaces via capacitary Lorentz spaces associated with Besov capacities.
For this purpose, we first study the embedding of the associated capacitary Lorentz space to the classical Lorentz space. Then,
 the embedding of the Besov space to the capacitary Lorentz space is established. Meanwhile, we show that
these embeddings are closely related to the iso-capacitary type inequalities in terms of a new-introduced
fractional $(\beta, p, q)$-perimeter. Moreover, characterizations of more general Sobolev type inequalities in Besov spaces have also been established.
\end{abstract}

\maketitle

\tableofcontents \pagenumbering{arabic}


\section{Introduction}
The Sobolev inequality plays a significant role in harmonic analysis, mathematical physics, and PDEs. Interested readers are referred to see \cite{Beckner, Carlen, Del Pino, Gross, Merker} and the references therein for more about Sobolev type inequalities.   Fractional calculus and fractional PDEs have become extremely  popular in mathematics, physics, engineering science and other areas, see \cite{Caffarelli,Frank} for instance. Fractional Sobolev inequalities have been studied in many references, see \cite{Cotsiolis, Cotsiolis2, Cotsiolis3, Hajaiej Yu Zhai, Ludwig, Palatucci,  Xiao 2016, Xiao20061,  Xiao2006, Xiao20062,   Xiao Zhai 2} for instance.

Denote by $C^{\infty}_{0}(\mathbb R^{n})$  the set of all compactly supported infinitely differentiable functions. Let $1\leq p<n/\beta$. The following
fractional Sobolev type inequality in Besov spaces $\dot{\Lambda}_\beta^{p,p}(\mathbb{R}^n)$ holds:
  \begin{eqnarray}\label{ClassSobolev}
\left(\int_{\mathbb{R}^n}|f(x)|^{\frac{np}{n-p\beta}}dx\right)^{\frac{n-p\beta}{np}}
\lesssim\|f\|_{\dot{\Lambda}_\beta^{p,p}(\mathbb{R}^n)}\quad
 \forall\ f\in C_0^\infty(\mathbb{R}^n),
 \end{eqnarray}
see \cite[Theorem 7.34]{AF} and \cite[Theorem 6.5]{NPV}. Motivated by \cite{Xiao2007, Xiao 2016}, in this paper, we will strengthen the inequality (\ref{ClassSobolev}) as follows:
 \begin{eqnarray}\label{Sobolev}
\left(\int_{\mathbb{R}^n)}|f(x)|^{\frac{np}{n-p\beta}}dx\right)^{\frac{n-p\beta}{np}}
&\lesssim&\left(\int_0^\infty V(O_t(f))^{\frac{n-p\beta}{n}}dt^p\right)^{1/p}\\
&\lesssim& \|f\|_{\dot{\Lambda}_\beta^{p,p}(\mathbb{R}^n)}\quad
 \forall\ f\in C_0^\infty(\mathbb{R}^n). \nonumber
 \end{eqnarray}
 Here $O_t(f)=\left\{x\in\mathbb{R}^n:\ |f(x)|>t\right\},$ and $V(O)$ denotes the volume of $O$ which is defined as the Lebesgue integral of the characteristic function $1_{O}$ on $\mathbb R^{n}$, i.e.,
 $$V(O):=\int_{\mathbb R^{n}}1_{O}(x)dx.$$

 On the other hand, we will show that   fractional Sobolev type inequalities in Besov spaces   can be strengthened by a Choquet  integral with respect to  Besov capacities.  That is, the second term in (\ref{Sobolev}) can be replaced by  capacitary Lorentz norms:
\begin{eqnarray}\label{stengthentwo}
\left(\int_{\mathbb{R}^n}|f(x)|^{\frac{np}{n-p\beta}}dx\right)^{\frac{n-p\beta}{np}}&\lesssim& \left(\int_{0}^\infty \left(C^{p,p}_{\beta}(O_t(f))\right)^{\frac{n}{n-p\beta}}dt^{\frac{np}{n-p\beta}}\right)^{\frac{n-p\beta}{np}}\\
&\lesssim& \|f\|_{\dot{\Lambda}_\beta^{p,p}(\mathbb{R}^n)}\quad \forall\  f\in C_0^\infty(\mathbb{R}^n).\nonumber
    \end{eqnarray}
    We  will also show that  (\ref{Sobolev}) is equivalent to the iso-capacitary type inequality in terms of Besov capacities.  While  (\ref{stengthentwo}) implies the iso-capacitary type inequalities in terms of Besov capacities and a new introduced fractional $(\beta,p,q)-$perimeter $P^{p,q}_\beta(\cdot).$

Moreover,  we will show that the general Sobolev type inequality 
\begin{equation}\label{geenralSobopq1}
    \|f\|_{L^{\frac{np}{n-p\beta},q}(\mathbb{R}^n)}\lesssim \|f\|_{\dot{\Lambda}^{p,q}_\beta(\mathbb{R}^n)},
    \end{equation} established in \cite[Theorem 7.34]{AF}, is equivalent to a new fractional Hardy inequality, an  iso-capacitary inequality,  and can be strengthened  by capacitary Lorentz norms when $q> p.$

Here, for $\beta\in (0,n),$  $p\in(0,n/\beta)$ and $q>0,$ $\dot{\Lambda}^{p,q}_\beta(\mathbb{R}^n)$ are defined as the closure of all $C^\infty_0$ functions $f$ with $\|f\|_{\dot{\Lambda}^{p,q}_\beta(\mathbb{R}^n)}<\infty.$
The norm $\|f\|_{\dot{\Lambda}^{p,q}_\beta(\mathbb{R}^n)}$ is defined  as follows.
$$\|f\|_{\dot{\Lambda}^{p,q}_\beta(\mathbb{R}^n)}=\left(\int_{\mathbb{R}^n}\|\Delta ^k_hf\|^q_{L^p}|h|^{-(n+\beta q)}dh\right)^{1/q}.$$
Here $k=1+[\beta], \beta=[\beta]+\{\beta\}$ with $[\beta]\in \mathbb{Z}_+,$ $\{\beta\}\in(0,1)$ and
\begin{equation*}
\triangle^k_{h}f(x)=
\left\{\begin{aligned}
&\triangle^1_h\triangle^{k-1}_h f(x),\quad k>1;\\
&f(x+h)-f(x), \quad k=1.
\end{aligned}\right.
\end{equation*}

 For a compact set $K\subset \mathbb{R}^{n}$,  the Besov capacity $C^{p,q}_\beta(K)$ is defined as
$$C_{\beta}^{p,q}(K):=\inf\left\{\|f\|^{p}_{\dot{\Lambda}^{p,q}_\beta(\mathbb{R}^n)}:\ f\in C_0^{\infty}(\mathbb{R}^n)\text{ and } f\geq 1_K\right\}$$
and for any set $E\subset\mathbb R^{n}$, one defines
$$C_{\beta}^{p,q}(E):=\inf_{\text{ open } O\supseteq E}\sup_{\text{ compact }K\subseteq O}\left\{C^{p,q}_{\beta}(K)\right\},$$
where $1_E$ denotes the characteristic function of $E$. The Besov capacity $C^{p,q}_\beta(\cdot)$  has been studied in \cite{Adams1, Adams Xiao} for $1<q< \infty, $ in \cite{Xiao Zhai 2} for $p=q\in (0,1).$

This paper is organized as follows. In Section \ref{sec-2}, we will show that Besov spaces can be embedded to  capacitary Lorentz spaces.  In Section \ref{sec 3}, we first prove that  fractional Sobolev type inequalities are equivalent to  fractional Hardy inequalities and   iso-capacitary inequalities. Then, we get strengthened Sobolev inequalities by   Lorentz norms.   In Section \ref{sec 4}, firstly, we will study  the embeddings of capacitary Lorentz spaces  to  the classical Lorentz space, and that of  Besov space    to capacitary Lorentz spaces. Finally, we strengthen fractional Sobolev type inequalities by  capacitary Lorentz norms.

{\it Some notations}:
\begin{itemize}
\item $U\approx V$ means that there is a constant $C>0$ such that $C^{-1}V\leq U\leq C V$. If $U\lesssim V$, then we write $U\lesssim V$. Similarly, we write $V\gtrsim U$ if $V\geq CU$.
\item  Let $k\in\mathbb{N}\cup \{\infty\}$. The symbol  ${C}^{k}(\mathbb{R}^{n})$ denotes the class of all functions
  $f:\ \mathbb{R}^{n}\rightarrow \mathbb{R}$ with $k$ continuous partial derivatives.
  For any subset $E\subseteq\mathbb R^{n}$, denote by $1_{E}$ the characteristic function of $E$.

\end{itemize}

\section{Capacitary Lorentz Spaces}\label{sec-2}
In this paper, we need to use  Lorentz/Lebesgue spaces associated with a nonnegative Radon measure $\mu,$ or the Besov capacity $C^{p,q}_\beta(\cdot)$,
or the Netrusov capacity. The Netrusov capacity $H^{\varepsilon}_{d,\theta}(\cdot)$ with $(\varepsilon,d,\theta)\in (0,\infty]\times (0,\infty)\times (0,\infty)$ is defined as, see \cite{Netrusov},
$$H^{(\varepsilon)}_{d,\theta}(K)=\inf\left(\sum_{i=0}^\infty(m_i2^{-id})^\theta\right)^{1/\theta},$$
where the infimum is taken over all countable coverings of $K\subset \mathbb{R}^n$ by balls whose radii $r_j$ do not exceed $\varepsilon,$ while $m_i$ is the number of balls from this covering whose radii $r_j$ belong to the interval $(2^{-i-1},2^{-i}],$ $i=0,1,2,\ldots.$ When $\theta=1,$ the Netrusov capacity $H^{(\infty)}_{d,1}(\cdot)$ is  the classical Hausdorff capacity $H^{(\infty)}_{d}(\cdot).$

 Denote by $\nu$ either   a nonnegative Radon measure $\mu$ on $\mathbb{R}^{n},$ or the Besov capacity $C^{p,q}_\beta(\cdot)$ with $0<p,q<\infty,$ or the Netrusov capacity $H^{(\varepsilon)}_{d,\theta}(\cdot).$ For $0<p_0,q_0<\infty,$
$L^{p_0,q_0}(\mathbb{R}^{n},\nu)$ and  $L^{p_0}(\mathbb{R}^{n},\nu)$  denote the Lorentz space and the Lebesgue space of all functions $g$ on $\mathbb{R}^{n}$, respectively, for which
$$\|g\|_{L^{p_0,q_0}(\mathbb{R}^{n};\nu)}=\left(\int_0^\infty\left(\nu(O_t(g))\right)^{q_0/p_0}d\lambda^{q_0}\right)^{1/q_0}<\infty$$
and
$$\|g\|_{L^{p_0}(\mathbb{R}^{n};\nu)}=\left(\int_{\mathbb{R}^{n}}|g(x)|^{p_0}d\nu\right)^{1/p_0}<\infty,$$
respectively. Moreover, we denote by $L^{p_0,\infty}(\mathbb{R}^{n};\nu)$ the set of all $\nu-$measurable functions $g(\cdot)$ on $\mathbb{R}^{n}$ with
$$\|g\|_{L^{p_0,\infty}(\mathbb{R}^{n};\nu)}=\sup_{s>0}s\left(\nu\left(O_s(g)\right)\right)^{1/p_0}<\infty.$$

The following result is standard. For the readers' convenience, we provide the proof here.
\begin{lemma}\label{Lemma21} For $0<q_0\leq r<\infty \ \&\ p_0>0$, there holds
\begin{equation}\label{inclusionofBesovLorentz}
L^{p_0,q_0}\left(\mathbb{R}^{n};\nu\right)
\hookrightarrow
L^{p_0,r}\left(\mathbb{R}^{n};\nu\right)
\hookrightarrow    L^{p_0,\infty}\left(\mathbb{R}^{n};\nu\right).
\end{equation}
\end{lemma}
\begin{proof}
Since  $\nu(O_t(f))$ is monotone decreasing in $t,$
we have   $$\frac{d}{dt}\left(\int_0^t(\nu(O_s(f)))^{r/p_0}ds^r\right)^{p_0/r}\geq p_0\nu(O_t(f))t^{p_0-1}$$
   and
   $$(s^{p_0}\nu(O_s(f)))^{r/p_0}\leq
   \left(p_0\int_0^\infty \nu(O_t(f))t^{p_0-1}dt\right)^{r/p_0}\leq\int_0^\infty(\nu(O_t(f)))^{r/p_0}dt^{r}\quad \forall\ s>0,$$
   which implies
   $\|f\|_{L^{p_0,\infty}\left(\mathbb{R}^{n};\nu\right)}\lesssim \|f\|_{L^{p_0,r}\left(\mathbb{R}^{n};\nu\right)}.
$
Thus, we have $
L^{p_0,r}\left(\mathbb{R}^{n};\nu\right)
\hookrightarrow    L^{p_0,\infty}\left(\mathbb{R}^{n};\nu\right).$

On the other hand, since  $$\|f\|_{L^{p_0,r}\left(\mathbb{R}^{n};\nu\right)}=\left(\int_0^\infty\left(\nu\left(\{O_\lambda(f)\right)\right)^{r/p_0}d\lambda^{r}\right)^{1/r},
$$
we have
\begin{eqnarray*}
\|f\|_{L^{p_0,r}\left(\mathbb{R}^{n};\nu\right)}&=&\left(\int_0^\infty\lambda^{r-q_0+q_0}
\left(\nu\left(O_\lambda(f)\right)\right)^{{r}/{p_0}-{q_0}/{p_0}+{q_0}/{p_0}}\frac{d\lambda}{\lambda}\right)^{1/r}\\
&\lesssim&\left(\sup_{\lambda>0}\left(\lambda \nu\left(O_\lambda(f)\right)\right)^{1/p_0}\right)^{1-{q_0}/{r}}\left(\int_0^\infty\lambda^{q_0}\left(\nu\left(O_\lambda(f)\right)\right)^{{q_0}/{p_0}}\frac{d\lambda}{\lambda}\right)^{1/r}\\
&\lesssim&\|f\|_{L^{p_0,\infty}\left(\mathbb{R}^{n};\nu\right)}^{1-{q_0}/{r}}\|f\|_{L^{p_0,q_0}\left(\mathbb{R}^{n};\nu\right)}^{q_0/r}\\
&\lesssim&\|f\|_{L^{p_0,r}\left(\mathbb{R}^{n};\nu\right)}^{1-{p_0}/{r}}\|f\|_{L^{p_0,q_0}\left(\mathbb{R}^{n};\nu\right)}^{q_0/r},
\end{eqnarray*}
which gives $\|f\|_{L^{p_0,r}\left(\mathbb{R}^{n};\nu\right)}\lesssim \|f\|_{L^{p_0,q_0}\left(\mathbb{R}^{n};\nu\right)}$ and thus $L^{p_0,q_0}\left(\mathbb{R}^{n};\nu\right)
\hookrightarrow
L^{p_0,r}\left(\mathbb{R}^{n};\nu\right).$
\end{proof}

When $\nu=C^{p,q}_\beta(\cdot)$, we call $L^{p_0,q_0}(\mathbb{R}^{n};C^{p,q}_\beta),$  $L^{p_0,\infty}(\mathbb{R}^{n};C^{p,q}_\beta),$ and $L^{p_0}(\mathbb{R}^{n};C^{p,q}_\beta)$     the associated capacitary  Lorentz/Lebesgue  space.
We can show that  Besov spaces are embedded to the associated capacitary Lorentz spaces.
Let $p\lor q=\max\{p,q\}$ and $p\land q=\min\{p,q\}.$

 \begin{proposition}\label{strongestimate}
    Let $\beta\in (0,1)$, $p=q\in\left({n}/{(n+\beta)},1\right),$ or $(\beta,p,q)\in(0,n)\times [1,n/\beta)\times (1,\infty),$ $p\lor q\leq r\leq \infty.$
There holds
    \begin{equation}\label{sci}
        \|f\|_{L^{p,r}\left(\mathbb{R}^{n};C^{p,q}_\beta\right)}
\lesssim \|f\|_{\dot{\Lambda}^{p,q}_\beta(\mathbb{R}^n)}
    \quad\forall\ f \in C_0^\infty(\mathbb{R}^n).
        \end{equation}
\end{proposition}

\begin{proof}
The case $r=p\lor q$
is due to Maz'ya \cite{Mazya} when $p=q>1.$ When $1\leq p\leq q<\infty, 0<\beta<1,$ Wu \cite{Wu} proved (\ref{sci}).  Adams-Xiao \cite{Adams Xiao} established (\ref{sci}) when $0<\beta<\infty, (p,q)\in (1, n/\beta)\times (1,\infty).$ Xiao-Zhai
\cite{Xiao Zhai 2} shown that (\ref{sci}) holds when $0<\beta<1, {n}/{(n+\beta)}<p=q<1.$
The case $r>p\lor q$ follows from  (i) and  the inclusion (\ref{inclusionofBesovLorentz}) in Lemma \ref{Lemma21} directly.
\end{proof}

Proposition \ref{strongestimate} is very important in studying Sobolev type inequalities and  Carleson embeddings problems, see \cite{Adams Xiao,  Li2020, Li, Xiao Zhai 2} for instance. We will use Proposition \ref{strongestimate} to establish  the main results in this paper. Moreover,  it  helps us to get the following result which generalizes Admas' inequalities  to general fractional Besov space.
In \cite{Adams86}, Adams proved that
$$\int_0^\infty H^{(\infty)}_{n-k}\left(O_t(f)\right)dt\lesssim \|\nabla^kf\|_{L^1(\mathbb{R}^n)}\quad \forall\ f\in C^\infty_0(\mathbb{R}^n),$$
which was generalized  by Xiao in \cite{Xiao20062} to  endpoint Besov spaces:
$$\int_0^\infty H^{(\infty)}_{n-\beta}\left(O_t(f)\right)dt\lesssim \|f\|_{\dot{\Lambda}^{1,1}_\beta(\mathbb{R}^n)}\quad \forall\ f\in C^\infty_0(\mathbb{R}^n),$$
when $\beta\in (0,n).$
Here, we will get a similar inequality for $\dot{\Lambda}^{p,q}_\beta(\mathbb{R}^n)$  and the Netrusov capacity.
\begin{proposition}\label{Pro1}
     Let  $\beta\in (0,1)$, $p=q\in\left({n}/{(n+\beta)},1\right),$ or $(p,q)\in \left(1,{n}/{\beta}\right)\times (1,\infty).$ If $r\in [p\lor q, \infty],$ then there holds $$
     \|f\|_{L^{p,r}\left(\mathbb{R}^{n};H^{(\infty)}_{n-p\beta,q/p}\right)}
  \lesssim \|f\|_{\dot{\Lambda}^{p,q}_\beta(\mathbb{R}^n)}
    \quad\forall\ f \in C_0^\infty(\mathbb{R}^n)
.$$
\end{proposition}
\begin{proof}
It follows from \cite[Theorem 2]{Adams2} and  \cite[Theorem 2]{Netrusov} that $C^{p,q}_\beta(\cdot)\approx H^{(\infty)}_{n-p\beta,q/p}(\cdot).$
Thus, (ii) of
 Proposition \ref{strongestimate} finishes the proof.
\end{proof}

\section{Strengthened  Fractional Sobolev Inequalities by Lorentz Spaces}\label{sec 3}
We will discuss the case $p=q$ and $p\neq q$ separately. 
\subsection{The case $p=q$}
We will show that the fractional Sobolev type inequality in $\dot{\Lambda}_\beta^{p,q}(\mathbb{R}^n)$ can be strengthened. Our first result holds for general $p$ and $q.$  
Recall $p\lor q=\max\{p,q\}$ and $p\land q=\min\{p,q\}.$

 \begin{theorem}\label{them2}
 Let
 $\beta\in (0,1)$ and  $p=q\in\left({n}/{(n+\beta)},1\right),$ or $(\beta,p,q)\in(0,n)\times [1, {n}/{\beta})\times (0,\infty).$ Let $p_0\geq q_0\geq p\lor q,$
and $\mu$ be a non-negative Radon measure.  Then the following statements are equivalent.
\item{\rm (i)} For any $ f\in C_0^\infty(\mathbb{R}^n)$,
$$\|f\|_{L^{p_0, q_0}(\mathbb{R}^n;\mu)}\lesssim \|f\|_{\dot{\Lambda}_\beta^{p,q}(\mathbb{R}^n)}.$$
\item{\rm (ii)} For any $ f\in C_0^\infty(\mathbb{R}^n)$,
$$\|f\|_{L^{p_0}(\mathbb{R}^n;\mu)}\lesssim\|f\|_{\dot{\Lambda}_\beta^{p,q}(\mathbb{R}^n)}.$$
\item{\rm (iii)} For any $ f\in C_0^\infty(\mathbb{R}^n)$,
$$\|f\|_{L^{p_0,\infty}(\mathbb{R}^n;\mu)}\lesssim\|f\|_{\dot{\Lambda}_\beta^{p,q}(\mathbb{R}^n)}.$$
\item{\rm(iv)} For any bounded domain $O\subset \mathbb{R}^n$ with $C^{\infty}$ boundary $\partial O$,
 $$\left(\mu(O)\right)^{p/p_0}\lesssim C^{p,q}_{\beta}(\overline{O}).$$

\end{theorem}
  \begin{proof}
Since  the Lorentz space $L^{r,q}(\mathbb R^{n};\mu)$ is increasing with respect to $q$ and
$p_0\geq q_0$, we have
 the implications (i)$\Longrightarrow$(ii)$\Longrightarrow$(iii).

(iii)$\Longrightarrow$(iv). Assume that (iii) is true. Given a   bounded domain $O\subset \mathbb{R}^n,$ for any non-negative $f\in C_0^\infty(\mathbb{R}^n) $ with $f\geq 1$ on ${O}$, i.e., ${O}\subset \overline{O_1(f)},$ we have
\begin{equation}\label{3.1new}
    \mu(O)\leq \mu(\overline{O_1(f)})\lesssim\|f\|^{p_0}_{\dot{\Lambda}_\beta^{p,q}(\mathbb{R}^n)},
    \end{equation}
 which implies    (iv) by taking infimum on $f$ on the right hand side of (\ref{3.1new}).

(iv)$\Longrightarrow$(i). Assume that (iv) holds. Then,  one has $(\mu(O))^{q_0/p_0}\lesssim (C^{p,q}_{\beta}(\overline{O}))^{q_0/p}$ and hence,
\begin{eqnarray}\label{eq-3.1}
\|f\|_{L^{p_0,q_0}(\mathbb{R}^n;\mu)}&=&\left(\int_0^\infty\left(\mu(O_t(f))\right)^{q_0/p_0}dt^{q_0}\right)^{1/q_0}\\
&\lesssim&\left(\int_0^\infty \left(C^{p,q}_\beta(O_t(f))\right)^{q_0/{p}}dt^{q_0}\right)^{1/q_0}\nonumber\\
&\lesssim&\|f\|_{\dot{\Lambda}_\beta^{p,q}(\mathbb{R}^n)},\nonumber
\end{eqnarray}
where in the last inequality of (\ref{eq-3.1}), we have used  Proposition \ref{strongestimate} since $q_0\geq p\lor q.$ This indicates (i).
\end{proof}

Theorem \ref{them2} itself is very important because it characterizes  a Radon measure such that the Sobolev type inequality holds.
The case of $p=q=q_0$ of Theorem \ref{them2} has been studied in  \cite{Xiao20062,Zhai,Xiao Zhai 2}.
  In  \cite{Frank Seiringer}, Frank and Seiringer proved that the fractional Sobolev inequality can be deduced from the fractional Hardy inequality, see \cite[Theorem 4.1]{Frank Seiringer}.  Using the case $p=q=q_0$ and $p_0={np}/{(n-\beta)}$ of Theorem \ref{them2}, we  prove   the following theorem which shows the equivalent of the fractional Sobolev inequality,   the fractional Hardy inequality, and the iso-capacitary inequality.
\begin{theorem}\label{ThemequalyofSobolevHardy}
 Let $\beta\in (0,1), 1\leq p<{n}/{\beta}.$  Then the following statements are equivalent.

\item{\rm (i)} The analytic inequality:
\begin{equation}\label{IneqstrongSolobevInequ}
\left(\int_0^\infty (V\left(O_t(f)\right))^{\frac{n-p\beta}{n}}dt^p\right)^{1/p}\lesssim
\|f\|_{\dot{\Lambda}^{p,p}_\beta(\mathbb{R}^n)}
    \quad\forall\ f \in C_0^\infty(\mathbb{R}^n).
\end{equation}

\item{\rm (ii)} The fractional Sobolev inequality:  \begin{equation}\label{f-SobolevIneq}
\left(\int_{\mathbb{R}^n}|f(x)|^{\frac{np}{n-p\beta}}dx\right)^{\frac{n-p\beta}{np}}\lesssim \|f\|_{\dot{\Lambda}_\beta^{p,p}(\mathbb{R}^n)}\quad
 \forall\ f\in C_0^\infty(\mathbb{R}^n). \end{equation}

\item{\rm (iii)} The fractional Hardy inequality:
\begin{equation}\label{fractionalHaydy}
    \left(\int_{\mathbb R^{n}}\frac{|f(x)|^p}{|x|^{p\beta}}dx\right)^{1/p}\lesssim \|f\|_{\dot{\Lambda}^{p,p}_\beta(\mathbb{R}^n)}
    \quad\forall\ f\in C^\infty_0(\mathbb{R}^n).
\end{equation}
    \item{\rm (iv)} The  iso-capacitary inequality: for any bounded domain $O\subset \mathbb{R}^n$ with $C^\infty$ boundary $\partial O$, 
\begin{equation}\label{iso-capacitary Inequalities}
(V(O))^{\frac{n-p\beta}{n}}\lesssim C^{p,p}_\beta(\overline{O}).
\end{equation}
Moreover, (\ref{IneqstrongSolobevInequ}), (\ref{f-SobolevIneq}), (\ref{fractionalHaydy}) and  (\ref{iso-capacitary Inequalities}) are all true.

\end{theorem}

\begin{proof}
In  Theorem \ref{them2}, let $q_0=p=q$, $p_0={np}/{n-p\beta}$ and $\mu$ be the Lebesgue measure on $\mathbb{R}^n$. We can see that  the equivalence of (i)$\Longleftrightarrow$(ii)$\Longleftrightarrow$(iv) holds. We will only provide a proof of (i)$\Longleftrightarrow$(iii).

For (i)$\Longrightarrow$(iii), we assume that (\ref{IneqstrongSolobevInequ}) holds.
If $A\subset \mathbb R^n$ is a Borel set of finite Lebesgue measure, we define $A^{\#}$, the
symmetric rearrangement of the set $A$, to be the open ball centered at
the origin whose volume is that of $A$. The symmetric-decreasing rearrangement, $f^{\#}$,
of a function $f$ is defined as follows. The symmetric-decreasing rearrangement of a characteristic function of
a set is obvious, namely, $1^{\#}_{A}=1_{A^{\#}}$. Now, if $f : \mathbb R^n \rightarrow\mathbb{C}$ is a Borel measurable function vanishing at infinity, we define
$$f^{\#}(x):=\int_{0}^{\infty}1^{\#}_{\{|f|>t\}}(x)dt$$
as the symmetric-decreasing rearrangement of $f.$

Firstly, we show that there exists a constant $C$ such that
\begin{equation}\label{specilEqu}
\left(\int_{\mathbb{R}^n}\frac{|f^{\#}(x)|^p}{|x|^{p\beta}}dx\right)^{1/p}=C\|f\|_{L^{\frac{np}{n-p\beta},p}},
\end{equation}
which is equivalent to  \cite[Lemma 4.3]{Frank Seiringer} since the Lorentz norm is invariant under the symmetric decreasing rearrangement.
For the convenience, we provide   the proof here. In fact, for $t>0$,
define $$O_t(f)=\Big\{x\in \mathbb{R}^n:\ |f(x)|>t\Big\}.$$ Then $V(O_t(f^p))=V(O_{t^{1/p}}(f))$. It follows from the equimeasurability of the functions $f$ and $f^{\#}$ that $V(O_{t}(f^{p}))=V(O_t((f^{\#})^p))$, which implies
\begin{eqnarray*}
\left(\int_{0}^\infty (V(O_t(f)))^{\frac{n-p\beta}{n}}dt^p\right)^{1/p}&=&\left(\int_0^\infty(V(O_t(f^p)))^{\frac{n-p\beta}{n}}dt\right)^{1/p}\\
&=&\left(\int_0^\infty(V(O_t((f^{\#})^p)))^{\frac{n-p\beta}{n}}dt\right)^{1/p}.
\end{eqnarray*}
Since $(f^{\#})^{p}$ is  a non-negative symmetric decreasing function, $|O_t((f^{\#})^p)|$ is the same as the volume of a ball $B_{r(t)}$ with radius $$r(t)=\left(\frac{V(|O_t((f^{\#})^p)|)}{\omega_n}\right)^{1/n}$$ with $\omega_n$ the surface area of the unit sphere $\mathbb{S}^{n-1}.$
Thus, there exists a constant $C$ such that
\begin{eqnarray}\label{SpecialEq2}
\left(\int_{\mathbb{R}^n}\frac{|f^{\#}(x)|^p}{|x|^{p\beta}}dx\right)^{1/p}&=&\left(\int_0^\infty\int_{\mathbb{R}^n}\frac{1_{B_{r(t)}}(x)}{|x|^{p\beta}}
dxdt\right)^{1/p}\\
&=&C \left(\int_0^\infty|O_t((f^{\#})^p)|^{\frac{n-p\beta}{n}} dt\right)^{1/p}\nonumber\\
&=&C\|f\|_{L^{\frac{np}{n-p\beta},p}(\mathbb{R}^n)}\nonumber
\end{eqnarray}
 since $\|\cdot\|_{L^{\frac{np}{n-p\beta},p}(\mathbb{R}^n)}$ is invariant under the symmetric decreasing rearrangement.

Note that \cite[Theorem 3.4]{LiebLoss} reads   as
$$\int_{\mathbb{R}^n} f(x)g(x)dx\leq \int_{\mathbb{R}^n}f^{\#}(x)g^{\#}(x)dx.$$
Notice that  $(\Phi\circ |f|)^{\#}=\Phi\circ (f)^{\#},$ where $\Phi(t)=t^p$ is non decreasing  for $t>0$. Since $|x|^{-p\beta}$ is symmetric-decreasing,  it follows from  (\ref{specilEqu})  and  (\ref{IneqstrongSolobevInequ})  that
$$\left(\int_{\mathbb{R}^n}\frac{|f(x)|^p}{|x|^{p\beta}}dx\right)^{1/p}\leq  \left(\int_{\mathbb{R}^n}\frac{|f^{\#}(x)|^p}{|x|^{p\beta}}dx\right)^{1/p}
\lesssim \|f\|_{L^{\frac{np}{n-p\beta},p}}
\lesssim \|f\|_{\dot{\Lambda}^{p,p}_\beta(\mathbb{R}^n)},$$
which gives us  (\ref{fractionalHaydy}). Thus (iii) is true.

For (iii)$\Longrightarrow$(i), we assume (\ref{fractionalHaydy}) holds. Since the Lorentz norm is invariant under the symmetric decreasing rearrangement, using (\ref{SpecialEq2}) one has
\begin{eqnarray*}
\left(\int_0^\infty V\left(O_t(f)\right)^{\frac{n-p\beta}{n}}dt^p\right)^{1/p}
&=&
\left(\int_0^\infty (V (O_t(f^{\#}) ))^{\frac{n-p\beta}{n}}dt^p\right)^{1/p}
\\
    &=&
    \left(\int_{\mathbb{R}^n}\frac{|f^{\#}(x)|^p}{|x|^{p\beta}}dx\right)^{1/p}\\
&\lesssim&\|f^{\#}\|_{\dot{\Lambda}^{p,p}_\beta(\mathbb{R}^n)}\\
&\lesssim&
\|f\|_{\dot{\Lambda}^{p,p}_\beta(\mathbb{R}^n)},
\end{eqnarray*}
where in the last step, we have used
 \cite[Theorem 9.2]{AlmgrenLieb} which states that the symmetric decreasing rearrangement is continuous under the fractional Besov norm  $\dot{\Lambda}^{p,p}_\beta(\mathbb{R}^n).$

   \end{proof}
\begin{remark}
    \item{(i)}
The sharp constant  of inequalities (\ref{f-SobolevIneq})  is only known when  $p=1$, see  \cite[Theoreom 4.1]{Frank Seiringer} or \cite[Theorem 4.10]{Brasco}.
\item{(ii)} The requirement $\beta\in (0,1)$ is only used when proving the equivalence of (i) and (iii). For   the equivalence of (i)$\Longleftrightarrow$(ii)$\Longleftrightarrow$(iv), $\beta$ can take values in $(0, n).$

\end{remark}

Since $p<\frac{np}{n-p\beta},$  we have
$$L^{\frac{np}{n-p\beta},p}(\mathbb{R}^n)\hookrightarrow L^{\frac{np}{n-p\beta}}(\mathbb{R}^n).$$ Therefore,  (\ref{IneqstrongSolobevInequ})  and (\ref{f-SobolevIneq})   imply  the following strengthened  fractional Sobolev type inequality in Besov spaces.

\begin{corollary} Let $\beta\in (0,n)$ and $1\leq p\leq n/\beta.$ For any $f\in C_0^\infty(\mathbb{R}^n)$, there holds
\begin{eqnarray}\label{strengthenedSobolevIneq}
\left(\int_{\mathbb{R}^n}|f(x)|^{\frac{np}{n-p\beta}}dx\right)^{\frac{n-p\beta}{np}}&\lesssim&\left(\int_0^\infty V\left(O_t(f)\right)^{\frac{n-p\beta}{n}}dt^p\right)^{1/p}\\
&\lesssim& \|f\|_{\dot{\Lambda}_\beta^{p,p}(\mathbb{R}^n)}.\nonumber
 \end{eqnarray}
 \end{corollary}

\subsection{The case $p\neq q$}
The inequality (\ref{strengthenedSobolevIneq}) strengthens the fractional Sobolev type inequality in Besove space  $\dot{\Lambda}_\beta^{p,q}(\mathbb{R}^n)$ when $p=q.$  In the following, we  consider the case $q\neq p.$ Firstly,  we need to establish the following result similar to  Theorem \ref{them2} without the requirement of $p_0\geq  p\lor q.$

 \begin{proposition}\label{prop5}
 Let
 $(\beta,p,q)\in(0,n)\times [1, {n}/{\beta})\times (1,\infty).$ Let $p_0>0$ and $q_0\geq p\lor q,$ $p_0>0$
and $\mu$ be a non-negative Radon measure.  Then the following statements are equivalent.
\item{\rm (i)} For any $ f\in C_0^\infty(\mathbb{R}^n)$,
$$\|f\|_{L^{p_0,q_0}(\mathbb{R}^n;\mu)}\lesssim \|f\|_{\dot{\Lambda}_\beta^{p,q}(\mathbb{R}^n)}.$$
\item{\rm (ii)} For any $ f\in C_0^\infty(\mathbb{R}^n)$,
$$\|f\|_{L^{p_0,\infty}(\mathbb{R}^n;\mu)}\lesssim\|f\|_{\dot{\Lambda}_\beta^{p,q}(\mathbb{R}^n)}.$$
\item{\rm(iii)} For any bounded domain $O\subset \mathbb{R}^n$ with $C^{\infty}$ boundary $\partial O$,
 $$\left(\mu(O)\right)^{p/p_0}\lesssim C^{p,q}_{\beta}(\overline{O}).$$
\end{proposition}
  \begin{proof}
Since  the Lorentz spaces $L^{p_0,q_0}(\mathbb R^{n};\mu)$ is increasing with respect to $q_0,$  we have
 the implications (i)$\Longrightarrow$(ii).

(ii)$\Longrightarrow$(iii). Assume that (ii) is true. Given a   bounded domain $O\subset \mathbb{R}^n,$ for any non-negative $f\in C_0^\infty(\mathbb{R}^n) $ with $f\geq 1$ on ${O}$, i.e., ${O}\subset \overline{O_1(f)},$ we have
$$\mu(O)\leq \mu(\overline{O_1(f)})\lesssim\|f\|^{p_0}_{\dot{\Lambda}_\beta^{p,q}(\mathbb{R}^n)},$$
 which implies    (iii) by taking infimum on $f$.

(iii)$\Longrightarrow$(i). Assume that (iv) holds. Then,  one has $(\mu(O))^{q_0/p_0}\lesssim (C^{p,q}_{\beta}(\overline{O}))^{q_0/p}$ and hence,
\begin{eqnarray}\label{eq-3.9}
\|f\|_{L^{p_0, q_0}(\mathbb{R}^n;\mu)}&=&\left(\int_0^\infty\left(\mu(O_t(f))\right)^{q_0/p_0}dt^{q_0}\right)^{1/q_0}\\
&\lesssim&\left(\int_0^\infty \left(C^{p,q}_\beta(O_t(f))\right)^{q_0/p}dt^{q_0}\right)^{1/q_0}\nonumber\\
&\lesssim&\|f\|_{\dot{\Lambda}_\beta^{p,q}(\mathbb{R}^n)},\nonumber
\end{eqnarray}
where in the last inequality of (\ref{eq-3.9}), we have used  Proposition \ref{strongestimate} since $q_0\geq p\lor q.$ This indicates (i).
\end{proof}

\begin{theorem}\label{Theme36}
 Let  $(\beta,p,q)\in(0,n)\times [1, {n}/{\beta})\times [1,\infty).$  Then the following statements are equivalent.

\item{\rm (i)} The analytic inequality:
\begin{equation}\label{IneqstrongSolobevInequpq}
\left(\int_0^\infty (V\left(O_t(f)\right))^{\frac{(n-p\beta)(p\lor q)}{np}}dt^{p\lor q}\right)^{\frac{1}{p\lor q}}\lesssim
\|f\|_{\dot{\Lambda}^{p,q}_\beta(\mathbb{R}^n)}
    \quad\forall\ f \in C_0^\infty(\mathbb{R}^n).
\end{equation}
   \item{\rm (ii)} The  iso-capacitary inequality: for any bounded domain $O\subset \mathbb{R}^n$ with $C^\infty$ boundary $\partial O$, 
\begin{equation}\label{iso-capacitary Inequalitiespq}
(V(O))^{\frac{n-p\beta}{n}}\lesssim C^{p,q}_\beta(\overline{O}).
\end{equation}
\item{\rm (iii)} The fractional Hardy inequality: for $\gamma=n(1-(p\lor q)/p)+\beta(p\lor q)$,
\begin{equation}\label{eq-3.1222}
\left(\int_{\mathbb R^{n}}\frac{|f(x)|^{p\lor q}}{|x|^{\gamma}}dx\right)^{1/(p\lor q)}\lesssim \|f\|_{\dot{\Lambda}^{p,q}_\beta(\mathbb{R}^n)}\quad\forall\ f \in C_0^\infty(\mathbb{R}^n).
\end{equation}
 Moreover, when $q>p,$ (\ref{IneqstrongSolobevInequpq}),  (\ref{iso-capacitary Inequalitiespq}) and (\ref{eq-3.1222}) are all true.
\end{theorem}
\begin{proof} The equivalence of
 (\ref{IneqstrongSolobevInequpq}) and  (\ref{iso-capacitary Inequalitiespq}) is a special of Proposition \ref{prop5} when $\mu$ is the Lebesgue measure, $q_0=p\lor q,$ and $p_0={np}/{n-p\beta}.$   It follows from  \cite[Theorem 7.34]{AF} that $\dot{\Lambda}^{p,q}_\beta(\mathbb{R}^n) \hookrightarrow L^{\frac{np}{n-p\beta},q}(\mathbb{R}^n)$ when $q\in[1,\infty).$ Thus,  (\ref{IneqstrongSolobevInequpq}), (\ref{iso-capacitary Inequalitiespq})    and (\ref{iso-capacitary Inequalitiespq}) are   true when $q>p.$

Below we prove (i) is equivalent to (iii). We  first assume that (i) holds. For any Borel set $A\subset\mathbb R^{n}$ with finite Lebesgue measure, denote by
$A^{\#}$ the symmetric rearrangment of $A$, which is the open ball centered at the origin whose volume is that of $A$.  Let  $f^{\#}$ denote the symmetric-decreasing rearrangement of a function $f$. For $t>0$, define
$O_{t}(f):=\Big\{x\in\mathbb R^{n}:\ |f(x)|>t\Big\}.$ It is easy to see that $V(O_{t}(f^{p}))=V(O_{t^{1/p}}(f)), \ p>1$. Then we can get
\begin{eqnarray*}
&&\left(\int^{\infty}_{0}\left(V(O_{t}(f))\right)^{\frac{(n-p\beta)(p\lor q)}{np}}dt^{p\lor q}\right)^{\frac{1}{p\lor q}}\\
&=&\left(\int^{\infty}_{0}\left(V(O_{t}(f^{\#}))\right)^{\frac{(n-p\beta)(p\lor q)}{np}}dt^{p\lor q}\right)^{\frac{1}{p\lor q}}\\
&\approx&\left(\int^{\infty}_{0}\left(V(O_{t}((f^{\#})^{p\lor q}))\right)^{\frac{(n-p\beta)(p\lor q)}{np}}dt\right)^{\frac{1}{p\lor q}}.
\end{eqnarray*}
Denote by $B_{r(t)}$ the ball centered at the origin with the radius
$$r(t):=\left(V(O_{t}((f^{\#})^{p\lor q}))\right)^{1/n}.$$
Then
\begin{eqnarray*}
\left(\int_{\mathbb R^{n}}\frac{|f^{\#}(x)|^{p\lor q}}{|x|^{\gamma}}dx\right)^{\frac{1}{p\lor q}}
&=&\left(\int^{\infty}_{0}\int_{\mathbb R^{n}}\frac{1_{B_{r(t)}}(x)}{|x|^{\gamma}}dxdt\right)^{\frac{1}{p\lor q}}\\
&= &\left(\int^{\infty}_{0}\left(V(O_{t}((f^{\#})^{p\lor q}))\right)^{\frac{n-\gamma}{n}}dt\right)^{\frac{1}{p\lor q}}.
\end{eqnarray*}
Using \cite[Theorem 3.4]{LiebLoss}, we obtain
\begin{eqnarray*}
\left(\int_{\mathbb R^{n}}\frac{|f(x)|^{p\lor q}}{|x|^{\gamma}}dx\right)^{1/(p\lor q)}&\lesssim&\left(\int_{\mathbb R^{n}}\frac{|f^{\#}(x)|^{p\lor q}}{|x|^{\gamma}}dx\right)^{\frac{1}{p\lor q}}\\
&\lesssim&\left(\int^{\infty}_{0}\left(V(O_{t}((f^{\#})^{p\lor q}))\right)^{\frac{(n-p\beta)(p\lor q)}{np}}dt\right)^{\frac{1}{p\lor q}}\\
&\approx&\left(\int^{\infty}_{0}\left(V(O_{t}(f^{p\lor q}))\right)^{\frac{(n-p\beta)(p\lor q)}{np}}dt\right)^{\frac{1}{p\lor q}}\\
&\approx&\left(\int^{\infty}_{0}\left(V(O_{t}(f))\right)^{\frac{(n-p\beta)(p\lor q)}{np}}dt^{p\lor q}\right)^{\frac{1}{p\lor q}}\\
&\lesssim&\|f\|_{\dot{\Lambda}^{p,q}_\beta(\mathbb{R}^n)},
\end{eqnarray*}
which gives (\ref{eq-3.1222}).

Conversely, if (iii) holds, then via a similar procedure, we can deduce that
\begin{eqnarray*}
&&\left(\int^{\infty}_{0}\left(V(O_{t}(f))\right)^{\frac{(n-p\beta)(p\lor q)}{np}}dt^{p\lor q}\right)^{\frac{1}{p\lor q}}\\
&&\quad\approx \left(\int^{\infty}_{0}\left(V(O_{t}(f^{\#}))\right)^{\frac{(n-p\beta)(p\lor q)}{np}}dt^{p\lor q}\right)^{\frac{1}{p\lor q}}\\
&&\quad\approx \left(\int_{\mathbb R^{n}}\frac{|f^{\#}(x)|^{p\lor q}}{|x|^{\gamma}}dx\right)^{\frac{1}{p\lor q}}\\
&&\quad\lesssim \|f^{\#}\|_{\dot{\Lambda}^{p,q}_\beta(\mathbb{R}^n)}\\
&&\quad\lesssim \|f\|_{\dot{\Lambda}^{p,q}_\beta(\mathbb{R}^n)},
\end{eqnarray*}
which proves (i).
\end{proof}

\begin{remark} When $q>p,$ 
Theorem \ref{Theme36} implies that the
Sobolev inequality \begin{equation}\label{geenralSobopq}
    \|f\|_{L^{\frac{np}{n-p\beta},q}(\mathbb{R}^n)}\lesssim \|f\|_{\dot{\Lambda}^{p,q}_\beta(\mathbb{R}^n)},
    \end{equation} established in \cite[Theorem 7.34]{AF}, is equivalent to the iso-capacitary inequality:  
$(V(O))^{\frac{n-p\beta}{n}}\lesssim C^{p,q}_\beta(\overline{O}),$
and the fractional Hardy inequality:
$$\left(\int_{\mathbb R^{n}}\frac{|f(x)|^{ q}}{|x|^{n(1-q/p)+q\beta}}dx\right)^{1/q}\lesssim \|f\|_{\dot{\Lambda}^{p,q}_\beta(\mathbb{R}^n)}
$$
 which is more general than the fractional Hardy inequality in $\dot{\Lambda}^{p,p}_\beta(\mathbb{R}^n).$ In the next section, we will strengthen (\ref{geenralSobopq}) by capacitary Lorentz norms.
\end{remark}

\section{Strengthened Fractional Sobolev Inequalities by Capacitary Lorentz Spaces} \label{sec 4}

 \subsection{Embeddings of  Capacitrary Lorentz Spaces to Lorentz Spaces}
 In this section, we prove that  the second term in
(\ref{strengthenedSobolevIneq}) will be  replaced by  capacitary Lorentz norms.

\begin{theorem}\label{them4}
 Let $\beta\in (0,n), p\geq 1,$ $q>0,$ $1<r<\infty.$ Let $p_0\geq 1, q_0>0,$ and  $\mu$ be a non-negative Radon measure.  Then the following statements are equivalent.
\item{\rm (i)} The embedding:
    \begin{equation}\label{them 4-(1)}
\|f\|_{L^{r, q_0}(\mathbb{R}^n)}
\lesssim  \|f\|_{L^{p_0,q_0}(\mathbb{R}^n;C^{p,q}_{\beta})}
\quad \forall\  f\in C_0^\infty(\mathbb{R}^n).
    \end{equation}
\item{\rm (ii)} The iso-capacitary inequality:
 \begin{equation}\label{them 4-(2)}
(\mu(O))^{p_0/r}\lesssim C^{p,q}_{\beta}(\overline{O}),
    \end{equation}
holds  for any bounded domain $O\subset \mathbb{R}^n$ with $ C^\infty$ boundary $\partial O.$
\end{theorem}

  \begin{proof}
       Suppose that (\ref{them 4-(2)}) is true. For any $f\in C_0^\infty(\mathbb{R}^n),$  (\ref{them 4-(2)}) implies
       $$(\mu(O_t(f)))^{p_0/r}\lesssim  C^{p,q}_{\beta}(\overline{O_t(f)})$$ 
       and so
    \begin{eqnarray*}
    \|f\|_{L^{r,q_0}(\mathbb{R}^n)}&=&\left(\int_{0}^\infty \mu(O_t(f))^{q_0/r}dt^{q_0}\right)^{1/q_0}\\
    &\lesssim&  \left(\int_{0}^\infty (C^{p,q}_{\beta}(\overline{O_t(f)}))^{q_0/p_0}dt^{q_0}\right)^{1/q_0}.
    \end{eqnarray*}
Thus, (\ref{them 4-(1)}) holds.

Now, assume that (\ref{them 4-(1)}) is true. For any bounded domain $O\subset \mathbb{R}^n$ with $C^\infty$ boundary $\partial O,$ denote by $\hbox{dist}(x,E)$ the Euclidean distance of a point $x$ to a set $E.$ For any $\varepsilon\in (0,1), $ define

$$f_\varepsilon(x)=\left\{\begin{aligned}
    &1-\varepsilon^{-1}\hbox{dist}(x,\overline{O}),& \quad  \hbox{dist}(x,\overline{O})<\varepsilon;\\
     &0,&  \quad \hbox{otherwise}.
    \end{aligned}\right.$$
Thus, $f_\varepsilon\in C^\infty(\mathbb{R}^n)$ and
$$O_t(f_\varepsilon)\subset U_1:=\{x\in\mathbb{R}^n: \hbox{dist}(x,\overline{O})<1\}$$
and so $\mu(O_t(f_\varepsilon))\leq \mu(U_1)<\infty.$
Then, the dominated convergence theorem implies 
\begin{eqnarray}\label{proofofthem 4-(1)}
\lim_{\varepsilon\rightarrow 0^+}\|f_\varepsilon\|_{L^{r,q_0}(\mathbb{R}^n;\mu)}
&=&
\lim_{\varepsilon\rightarrow 0^+}\left(\int_0^1(\mu(O_t(f_\varepsilon)))^{q_0/p_0}dt^{q_0}\right)^{1/q_0}\nonumber\\
&=&
\left(\int_0^1\lim_{\varepsilon\rightarrow 0^+}(\mu(O_t(f_\varepsilon)))^{q_0/r}dt^{q_0}\right)^{1/q_0}\nonumber\\
&=&\left(\mu(O)\right)^{1/r}.
\end{eqnarray} On the other hand, (\ref{them 4-(1)}) implies
\begin{eqnarray*}
\|f_\varepsilon\|_{L^{r,q_0}(\mathbb{R}^n; \mu)}
&=&  \left(\int_{0}^\infty \left(\mu\left(O_t( f_\varepsilon)\right)\right)^{q_0/r}dt^{q_0}\right)^{1/q_0}\\
&\lesssim&  \left(\int_{0}^\infty \left(C^{p,q}_{\beta}\left(O_t( f_\varepsilon)\right)\right)^{q_0/p_0}dt^{q_0}\right)^{1/q_0}\\
&\lesssim&  \left(\int_{0}^1 \left(C^{p,q}_{\beta}\left(O_t( f_\varepsilon)\right)\right)^{q_0/p_0}dt^{q_0}\right)^{1/q_0}\\
&\lesssim & \left(C^{p,q}_{\beta}(\overline{O_\varepsilon})\right)^{1/p_0},
\end{eqnarray*}
where the last inequality is due to
$$\Big\{x\in\mathbb{R}^n:\ |f_\varepsilon(x)|\geq t\Big\}\subset \overline{O_\varepsilon}:=\Big\{x\in\mathbb{R}^n:\ \hbox{dist}(x,\overline{O}<\varepsilon)\Big\}.$$
 Letting $\varepsilon\rightarrow 0^+$ gives us
 (\ref{them 4-(2)}) due to  (\ref{proofofthem 4-(1)}).
\end{proof}

 \begin{theorem}\label{them3}
 Let $\beta\in (0,n), q>0,$  $1\leq p<{n}/{\beta},$ and $ q_0>0.$  Then the following two statements are equivalent.
\item{\rm (i)} The embedding:
    \begin{equation}\label{them 2-(1)}
\|f\|_{L^{\frac{np}{n-p\beta},q_0}(\mathbb{R}^n)}
\lesssim \|f\|_{L^{p,q_0}(\mathbb{R}^n; C^{p,q}_{\beta})}\quad \forall\  f\in C_0^\infty(\mathbb{R}^n).
    \end{equation}
\item{\rm (ii)} The iso-capacitary inequality
    \begin{equation}\label{them 2-(2)}
(V(O))^{\frac{n-p\beta}{n}}\lesssim  C^{p,q}_{\beta}(\overline{O})\quad \forall\ \hbox{bounded domain}\ O\subset \mathbb{R}^n\ \hbox{with}\ C^\infty\ \hbox{boundary}\ \partial O.
    \end{equation}
    Moreover, when $q\geq p,$ both  (\ref{them 2-(1)}) and (\ref{them 2-(2)}) are true.

 \end{theorem}

\begin{proof}
The truth of the iso-capacitary inequality   (\ref{them 2-(2)}) was established in
Theorem \ref{ThemequalyofSobolevHardy} for $p=q,$ and Theorem \ref{Theme36} for $q>p.$ So,
  (\ref{them 2-(2)}) implies the truth of  (\ref{them 2-(1)}) if we establish the equivalence of (\ref{them 2-(1)})-(\ref{them 2-(2)}) which is the special case of  Theorem \ref{them4} when $r={np}/({n-p\beta}),$ $p_0=p$ and $\mu$ is the Lebesgue measure on $\mathbb{R}^n.$

\end{proof}

 \subsection{Embeddings of Besov Spaces to Capacitrary Lorentz Spaces}
 
 Below we will show that  the embedding $\dot{\Lambda}^{p,q}_{\beta}(\mathbb{R}^n)\hookrightarrow L^{p,q_0}(\mathbb{R}^n, C^{p,p}_\beta)$  implies the iso-capacitary inequality in term of   a new introduced  fractional $(\beta,p,q)-$perimeter $P^{p,q}_\beta(E):$ 
 $$ \left(C^{p,q}_{\beta}(\overline{O})\right)^{1/p}\leq 2P^{p,q}_{\beta}(O)$$
 for $p\in [1,n/\beta),$ $\beta\in (0,1)$ and  all bounded domain $ O\subset \mathbb{R}^n$ with $C^\infty$boundary $\partial O$. Here $P^{p,q}_{\beta}(O)$  is defined as follows.
\begin{definition}\label{fractionalPerimeter} Let $p,q>0.$ For any  $E\subset \mathbb{R}^n,$ let $E^c=\mathbb{R}^n\backslash E$.
The fractional $(\beta,p,q)-$perimeter $P^{p,q}_\beta(E)$ is defined as
\begin{equation}\label{pqperimter}
P^{p,q}_\beta(E)=\left(\int_E\left(\int_{E^c}\frac{dx}{|x-y|^{\frac{(n+p\beta)p}{q}}}\right)^{q/p}dy\right)^{1/q}.
\end{equation}
\end{definition}
When $p=q=1, $ $P^{1,1}_\beta(E)=P_\beta(E)$ which is the fractional $\beta-$perimeter $P_\beta(E)$  defined as
 $$P_\beta(E):=\frac{1}{2}\|1_E\|_{\dot{\Lambda}^{1,1}_\beta(\mathbb{R}^n)}=\int_E\int_{E^c}\frac{dxdy}{|x-y|^{n+\beta}}.$$
The regularity   of set with minimal fractional perimeter $P_\beta(E),$ the approximation of $P_\beta(E)$ to the classical perimeter and other geometric properties of $P_\beta(E)$   have been studied in \cite{Brasco, CaffarelliRoquejoffre,Ambrosio, Bourgain, CaffarelliEnrico, Davila,PonceSpector, Ponce,Frank Seiringer,FuscoVincent}.
The fractional perimeter $P_\beta(E)$ has been applied to study other embeddings in \cite{Li}.
 When $p=q>1$ and $p\beta\in (0,1),$ $P^{p,p}_\beta(E)=(P_{p\beta}(E))^{1/p}.$

 \begin{theorem}\label{them5}
Let $\beta\in (0,n), 1\leq p<{n}/{\beta},$ $q>0,$ and $ p\lor q\leq q_0<\infty$.
    \item{\rm (i)} The following embedding holds:
        \begin{equation}\label{them 5-(1)}
\|f\|_{L^{p,q_0}(\mathbb{R}^n; C^{p,q}_{\beta})}\lesssim \|f\|_{\dot{\Lambda}_\beta^{p,q}(\mathbb{R}^n)}\quad \forall\  f\in C_0^\infty(\mathbb{R}^n).
    \end{equation}
\item{\rm (ii)} If $\beta\in(0,1),$ then the iso-capacitary inequality holds:  for any bounded domain $O\subset \mathbb{R}^n$ with $C^\infty$-boundary $\partial O$, there holds the following geometric inequality:    \begin{equation}\label{them 5-(2)}
 (C^{p,q}_{\beta}(\overline{O}))^{{1}/{p}}\lesssim P^{p,q}_{\beta}(O).
    \end{equation}

 \end{theorem}

\begin{proof}
For (i),    Proposition \ref{strongestimate} implies the truth of (\ref{them 5-(1)}). Now we prove (ii).  For  $\varepsilon>0$ and
 a bounded domain $O\subset \mathbb{R}^n$ with $C^\infty$ boundary $\partial O,$ denote $$O_\varepsilon:=\Big\{x\in\mathbb{R}^n:\ \text{ dist }(x,\overline{O}<\varepsilon)\Big\}$$
and
 $$f_\varepsilon(x):=\left\{\begin{aligned}
    &1-\varepsilon^{-1}\hbox{dist}(x,\overline{O}),&  \hbox{dist}(x,\overline{O})<\varepsilon;\\
     &0,& \hbox{otherwise}.
    \end{aligned}\right.$$
    Then, $f_\varepsilon(x)=1$ for all $x\in \overline{O}$ and so $\overline{O}\subset \overline{O_t(f_\varepsilon)}$ for all $\varepsilon\in(0,1)$ and $t\in (0,1).$ Thus,
    \begin{eqnarray*}
    (C^{p,q}_\beta(\overline{O}))^{1/p}\leq \left(\int_{0}^\infty \left(C^{p,q}_{\beta}(\overline{O_t(f_\varepsilon)})\right)^{q_0/p}dt^{q_0}\right)^{1/q_0}\lesssim \|f_\varepsilon\|_{\dot{\Lambda}_\beta^{p,q}(\mathbb{R}^n)}.
    \end{eqnarray*}
    Since $f_\varepsilon\rightarrow 1_{\overline{O}}$ as $\varepsilon\rightarrow 0^+,$ the dominated convergent theorem implies
    \begin{eqnarray*}
(C^{p,q}_\beta(\overline{O}))^{1/p}&\lesssim&\lim_{\varepsilon\rightarrow 0^+}\|f_\varepsilon\|_{\dot{\Lambda}_\beta^{p,q}(\mathbb{R}^n)}\\
&=&\|1_{\overline{O}}\|_{\dot{\Lambda}_\beta^{p,q}(\mathbb{R}^n)}\\
&=&\left(\int_{\mathbb{R}^n}\left(\int_{\mathbb{R}^n}\frac{|1_{\overline{O}}(x)-1_{\overline{O}}(y)|^p}{|x-y|^{(n+p\beta)p/q}}dx\right)^{q/p}dy\right)^{1/q}\\
&=&2\left(\int_{O^c}\left(\int_{O}\frac{1}{|x-y|^{(n+p\beta)p/q}}dx\right)^{q/p}dy\right)^{1/q}\\
&=&2P^{p,q}_{\beta}(O).
\end{eqnarray*}
Thus, (\ref{them 5-(2)}) holds.

 \end{proof}

 \begin{remark}
      When $p=q=1$ and $ q_0=n/(n-\beta),$ Xiao \cite{Xiao 2016} showed that (\ref{them 5-(1)}) and (\ref{them 5-(2)}) are true, sharp and  equivalent using the general co-area formula $$\|f\|_{\dot{\Lambda}_\beta^{1,1}(\mathbb{R}^n)}=2\int_{0}^\infty P_\beta(O_t(f))dt.$$

 \end{remark}

\subsection{Strengthened Fractional Sobolev  Inequalities by  Capacitary Lorentz Spaces}
Based on Theorems \ref{them3} and  \ref{them5}, we can strengthen the fractional Sobolev inequality $$\left(\int_{\mathbb{R}^n}|f(x)|^{\frac{np}{n-p\beta}}dx\right)^{\frac{n-p\beta}{np}}\lesssim  \|f\|_{\dot{\Lambda}_\beta^{p,p}(\mathbb{R}^n)}$$
and the isoperimetric type inequality
$$(V(O))^{1-p\beta/{n}}\lesssim   (2P^{p,p}_{\beta}(O))^p.$$

\begin{corollary}
Let $\beta\in (0,n)$ and $ 1\leq p<{n}/{\beta}, $ $p\lor q\leq q_0<\infty.$
\item{\rm (i)}  There holds  the analytic inequality: 
    \begin{equation}
    \label{them 6-(1)}
\|f\|_{L^{\frac{np}{n-p\beta},q_0}(\mathbb{R}^n)}\lesssim \|f\|_{L^{p,q_0}(\mathbb{R}^n; C^{p,q}_\beta)}\lesssim \|f\|_{\dot{\Lambda}_\beta^{p,q}(\mathbb{R}^n)}\quad \forall f\in C_0^\infty(\mathbb{R}^n).
   \end{equation}

      When $\beta\in (0,1),$  (i) implies the following geometric inequality.
\item{\rm (ii)} The geometric inequality:
    \begin{eqnarray*}\label{them 6-(2)}
(V(O))^{\frac{n-p\beta}{n}}\lesssim  C^{p,p}_{\beta}(\overline{O})\lesssim (P^{p,q}_{\beta}(O))^p\quad \forall\ \hbox{bounded domain}\ O\subset \mathbb{R}^n\ \hbox{with}\ C^\infty\ \hbox{boundary}\ \partial O.
    \end{eqnarray*}
 Thus, when  $\beta\in (0,1),$  (ii) is also true.
\end{corollary}
\begin{remark}

 \item{\rm (i)} When $q_0={np}/(n-p\beta)$ and $p=q,$ (\ref{them 6-(1)}) implies
    \begin{equation}\label{4.10}
 \left(  \int_{\mathbb{R}^n}|f(x)|^{\frac{np}{n-p\beta}}dx\right)^{\frac{n-p\beta}{np}}\lesssim \left(\int_{0}^\infty \left(C^{p,p}_{\beta}\left(O_t(f)\right)\right)^{\frac{n}{n-p\beta}}dt^{\frac{np}{n-p\beta}}\right)^{\frac{n-p\beta}{np}}\\
\lesssim \|f\|_{\dot{\Lambda}_\beta^{p,p}(\mathbb{R}^n)}
    \end{equation}
which strengthens  the fractional Sobolev inequality (\ref{ClassSobolev}
).  When $p=q=1,$ (\ref{4.10}) was established by Xiao in \cite{Xiao2007, Xiao 2016}.
\item{\rm (ii)} When  $q_0=q\geq p\geq 1,$
(\ref{them 6-(1)}) implies
$$\|f\|_{L^{\frac{np}{n-p\beta},q}(\mathbb{R}^n)}\lesssim \|f\|_{L^{p,q}(\mathbb{R}^n; C^{p,q}_\beta)}\lesssim \|f\|_{\dot{\Lambda}_\beta^{p,q}(\mathbb{R}^n)}
$$
which strengthens the Sobolev type inequality 
$$\|f\|_{L^{\frac{np}{n-p\beta},q}(\mathbb{R}^n)}\lesssim \|f\|_{\dot{\Lambda}_\beta^{p,q}(\mathbb{R}^n)}
$$
 established in \cite[Theorem 7.34]{AF}.

\item{\rm (iii)} When $p=1,$  in \cite{Hurri-Syrjanen},  $(V(O))^{\frac{n-\beta}{n}}\leq 2P^{1,1}_{\beta}(O)$ was proved to be equivalent to the fractional Sobolev inequality
$$\left(\int_{\mathbb{R}^n}|f(x)|^{\frac{np}{n-p\beta}}dx\right)^{\frac{n-p\beta}{np}}
\lesssim \|f\|_{\dot{\Lambda}_\beta^{p,p}(\mathbb{R}^n)}.
$$
\end{remark}

\end{document}